\documentclass[12pt]{article}
\pdfoutput=1
\usepackage{graphicx}
\usepackage{wrapfig}

\usepackage{amsmath,amsfonts,amssymb,amsthm}
\usepackage{txfonts}%
\usepackage{bm}
\usepackage{comment}

\theoremstyle{plain}

\theoremstyle{definition}

\theoremstyle{remark}

\begin{document}

\title{A power-law decay model with autocorrelation for posting data to social networking services
}

\author{
Toshifumi Fujiyama\thanks{Department of Mathematical Informatics, Graduate School of Information Science and Technology, University of Tokyo}, \ 
Chihiro Matsui\footnotemark[1]\ 
\ and Akimichi Takemura\footnotemark[1]
}
\date{November, 2014}
\maketitle

\begin{abstract}
We propose a power-law decay model with autocorrelation for posting data to social networking services
concerning particular events such as national holidays or major sport events.  
In these kinds of events we observe people's interest both before and after the events.
In our model the number of postings has a Poisson distribution whose expected value decays as a
power law.  Our model also incorporates autocorrelations by autoregressive specification of the
expected value.  We show that our proposed model well fits the data from social 
networking services\footnote{We are grateful to NTTCom Online Marketing Solutions Corporation for allowing us to use their BuzzFinder service, which provided excellent access to social networking service posting data.}.
\end{abstract}



\section{Introduction}
\label{sec:intro}

With the increasing usage of social networking services (SNS), such as many blog services, Facebook or Twitter,
it is important to obtain information from postings to these services. 
We look at the data on the number of postings for various events, such as national holidays 
or major sport events.
By analyzing the data, we gain insights on how a particular topic interests people, 
how actively it is discussed around the date of the event and how long the memory of the
event lasts.  Our data is time series data, but the length of the series is usually short,
covering several weeks before and after the event.  It is far from stationary, because
there is a sharp peak of the number of postings on the date of the event. 

Many studies on posting data to SNS focus on the network structure of the 
services (e.g.\ \cite{1}, \cite{3}, \cite{traud}), 
such as how certain topic spreads by interaction among users of the services.
In this paper we study the pattern of number of postings concerning particular events,
such as national holidays or major sport events.  The events we study are scheduled events, so
that people are aware when the events will happen.  This contrasts with unpredictable events,
such as the occurrence of large earthquakes (\cite{4}).  
In the case of unpredictable events, we only observe effects of the event after it happens.
On the other hand, in the case of the scheduled events, 
the number of postings depends on people's anticipation of the event
and the after-effect of the actual outcome of the event.
Another similar type of data studied in literature is the registration data until the deadline
for events such as academic conferences (e.g.\ \cite{alfi-et-al}, \cite{fenner}).
In this kind of data, we only observe people's actions before the deadline.

We combine the model of the mean number of postings proposed in \cite{5} 
and the conditional Poisson autoregressive (AR) models for count data.
The conditional Poisson autoregressive models are applied in many problems
in econometrics, political science or epidemiology (e.g.\ \cite{6}, \cite{7}).
For a theoretical survey of the conditional Poisson autoregressive model
see \cite{fokianos-JASA}, \cite{fokianos-handbook} and \cite{fokianos-2011}.
Zhu (\cite{zhu}) generalized the conditional Poisson model to  generalized Poisson integer-valued GARCH models.
By our model we can better predict how people's interest on particular events decays with time and 
our model is useful for example in designing advertisement strategy for web marketing.

\begin{figure}[htbp]
\centering
\includegraphics[width=6cm]{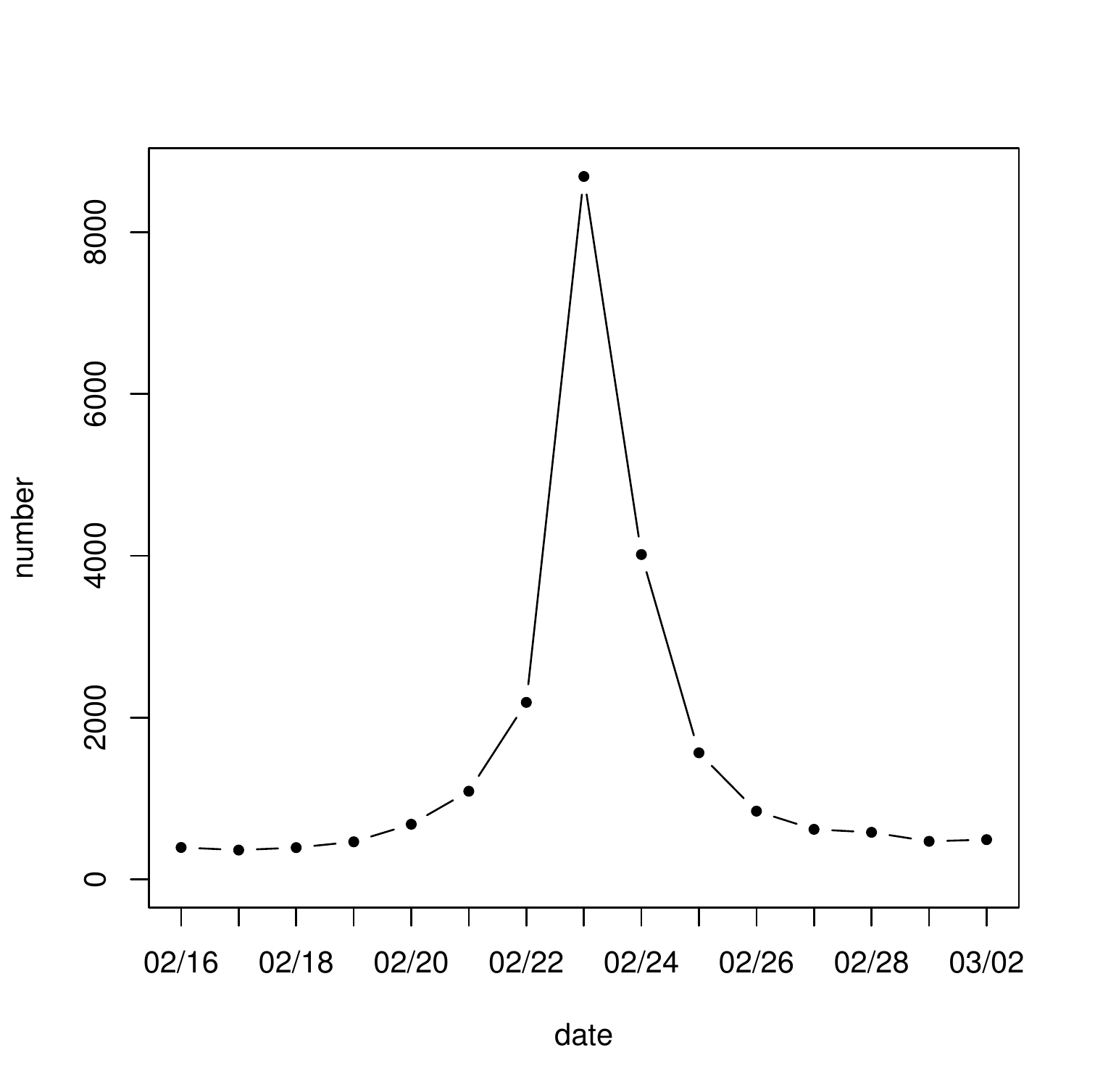} \quad
\includegraphics[width=6cm]{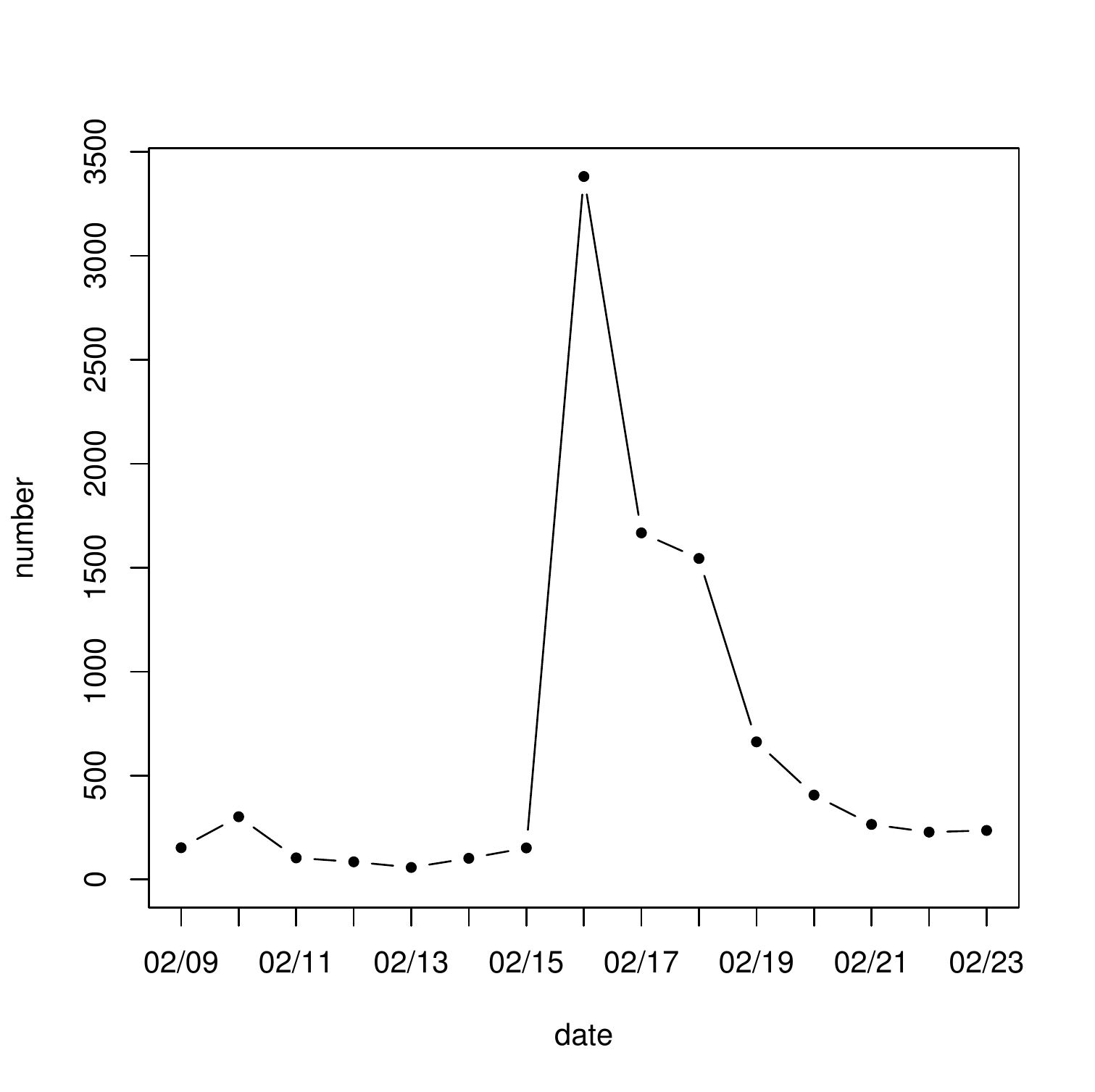}\\
Tokyo Marathon 2014 \hspace*{2.5cm} Silver medal by Kasai 
\caption{Symmetric and asymmetric patterns}
\label{fig:two-patterns}
\end{figure}

In the left graph of Figure \ref{fig:two-patterns} we show a typical symmetric pattern of number of postings.
Tokyo Marathon 2014 was held on Sunday, February 23 of 2014. It was well anticipated and it ended
without unexpected happenings. In these cases, the number of postings shows a symmetric
pattern with a sharp peak on the date of the event.  Both sides of the peak seems to
exhibit a behavior of some negative power in the time difference from the actual date of the event.
Also there is some stochastic fluctuation of the number of postings.
This is the motivation of our proposed model.
A different asymmetric pattern is observed, when there is some surprising element of the event.
In the right graph of Figure \ref{fig:two-patterns} we show the data on Noriaki Kasai around February 16, 2014, when he won a silver medal 
in ski jump in 2014 winter Olympics.  From the data we see that people did not anticipate the medal before
the event.

The organization of this paper is as follows.
In Section \ref{sec:power-decay-model} we propose to model the mean 
of the number of postings by a power-law decay model with three parameters and we assume
independent Poisson distributions.
We also study the Fisher information matrix under the proposed model.
In Section \ref{sec:conditional-Poisson} we incorporate autocorrelation  to our model
by conditional Poisson autoregression.  
In Section \ref{sec:data-analysis} we apply our model to some social networking service
data in Japan.  
We end the paper with some discussion in Section \ref{sec:summary}.

\section{A power-law decay model for the mean number of postings}
\label{sec:power-decay-model}
Let $t_0$ denote the date of the event and let $y_t$ denote the number of postings 
on the event on day $t$.
We model the expected value of $y_t$ by the following power-law decay function.
\begin{equation}
E(y_t) = \mu_t(\alpha,\beta,\gamma)=\gamma \frac{1}{(\alpha |t-t_0|+1)^\beta}.
\label{eq:power-decay}
\end{equation}
This power-law decay model was proposed by \cite{5} without the parameter $\alpha$.
They used the least-squares method, while we use the maximum likelihood estimation. 
They do not consider fitting the data close to the peak, which becomes possible by introducing the parameter $\alpha$.

The interpretation of the parameters is as follows.
\begin{description}
\setlength{\itemsep}{1pt}
\item $\alpha$: steepness of the curve just before and after the event
\item $\beta$: longer decay pattern
\item $\gamma$: impact of the event (peak level, the maximum number of postings)
\end{description}
In this section we assume that $y_t$, $t=t_L, t_L+1, \dots, t_U-1, t_U$, $t_L \le t_0 \le t_U$, 
are independent Poisson random variables with the mean given by  $\mu_t(\alpha,\beta,\gamma)$
in \eqref{eq:power-decay}. 
We call the model ``power-law decay independence model''.
We denote the probability function of the
Poisson distribution with the mean $\mu$ as
\begin{equation}
  \label{eq:poisson-probability-function}
     {\rm Po}(y \mid \mu)= \frac{\mu^y}{y!} e^{-\mu}.
\end{equation}
Then the likelihood function is written as
\begin{equation}
  L(\alpha,\beta,\gamma)
  =\prod_{t=t_L}^{t_U}{\rm Po}(y_{t} \mid \mu_t(\alpha,\beta,\gamma))
=
\prod_{t=t_L}^{t_U}
\frac{\mu_t(\alpha,\beta,\gamma)^{y_{t}}}{{y_{t}}!}
e^{-\mu_t(\alpha,\beta,\gamma)}.
\label{eq:independent-poisson-likelihood}
\end{equation}
\begin{figure}[htbp]
\centering
\includegraphics[width=7cm]{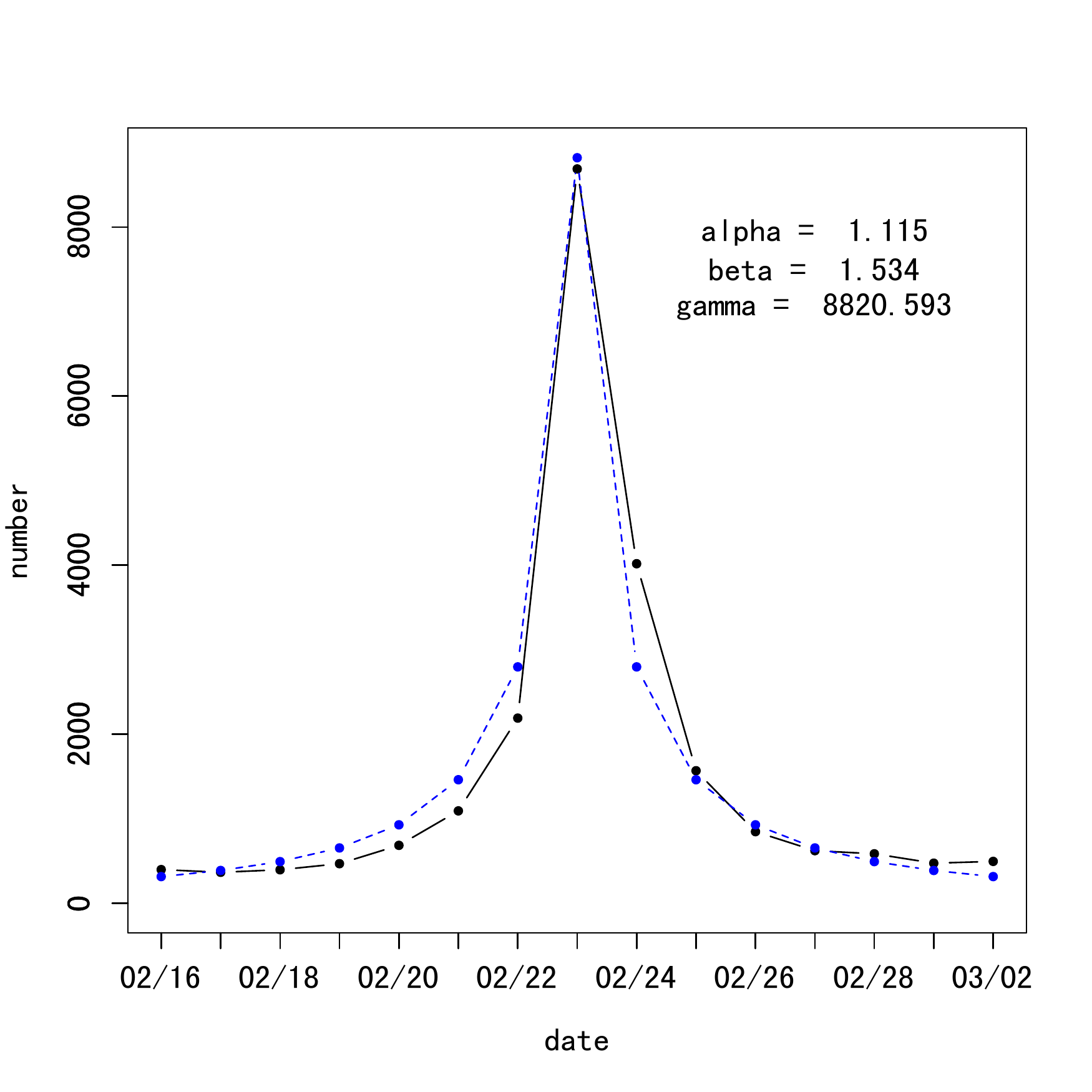} 
\caption{Fitting the power-law decay independence model to Tokyo Marathon data}
\label{fig:sample1}
\end{figure}

We found that the maximization of the log-likelihood function is numerically very simple.
We  give  a typical example of maximum likelihood estimation (MLE) of 
\eqref{eq:independent-poisson-likelihood} in Figure \ref{fig:sample1}.
In Figure \ref{fig:sample1} the solid line is the observed data and the 
the dotted line is the fitted curve of expected values by MLE.  We use the same distinction
of line types in Figures 3 and 4.
The estimates are $\hat\alpha=1.115$, $\hat\beta = 1.534$ and
$\hat\gamma=8820.593$.
Our model seems to fit the data well, but 
there is a slight asymmetry in this data, which is not captured by the symmetric model
in \eqref{eq:power-decay}.

We will discuss many examples of fitting 
of \eqref{eq:independent-poisson-likelihood} and its generalizations 
later in Section \ref{sec:data-analysis}.

\subsection{Fisher information matrix for the independence model} 
As we already noted in the beginning, the length of our data
is not very large and our data is far from stationary. 
Hence the usual asymptotics for the length of the series is not appropriate.
Nevertheless it is of theoretical interest to consider the behavior of MLE
of \eqref{eq:independent-poisson-likelihood} as 
the length of the series diverges to infinity.  
On the other hand we have large number of postings $y_{t_0}$ 
on the date $t_0$ of the event.  Under our model of Poisson distribution, 
we can also think of the asymptotics, where $\gamma=E(y_{t_0})$ diverges to infinity.
In the following we calculate the Fisher information matrix 
\[
I=
\begin{pmatrix}
I_{\alpha\alpha} & I_{\alpha\beta} & I_{\alpha\gamma} \\
I_{\alpha\beta} & I_{\beta\beta} & I_{\beta\gamma} \\
I_{\alpha\gamma} & I_{\beta\gamma} & I_{\gamma\gamma} 
\end{pmatrix},
\]
for our model to gain insights on the behavior 
of MLE for the model \eqref{eq:independent-poisson-likelihood}.

For notational simplicity let $t_0=0$ and $t_L \le 0 \le t_U$. Then
the log-likelihood function $l(\alpha,\beta,\gamma)=\log L(\alpha,\beta,\gamma)$
is written as
\begin{align*}
l(\alpha,\beta,\gamma) &= \sum_{t=t_L}^{t_U}  (y_t \log \mu_t(\alpha,\beta,\gamma) 
- \mu_t(\alpha,\beta,\gamma) - \log y_t!) \\
&= \sum_{t=t_L}^{t_U}  \Big(y_t (\log \gamma - \beta\log(\alpha |t|+1)) 
- \gamma  \frac{1}{(\alpha|t|+1)^\beta} - \log y_t!\Big) \\
&= C  + \log\gamma \sum_{t=t_L}^{t_U}  y_t -  \beta \sum_{t=t_L}^{t_U} y_t \log(\alpha |t| + 1)
 - \gamma \sum_{t=t_L}^{t_U} \frac{1}{(\alpha |t| +1)^\beta},
\end{align*}
where $C$ does not depend on parameters. Then
\[
-\frac{\partial^2}{\partial \alpha^2}l=
-\beta \sum_{t=t_L}^{t_U} y_t \frac{|t|^2}{(\alpha |t| + 1)^2}
+ \gamma \beta(\beta+1) \sum_{t=t_L}^{t_U} \frac{|t|^2}{(\alpha |t| +1)^{\beta+2}}.
\]
The expected value of this second derivative is given as
\begin{align*}
I_{\alpha\alpha}=E[-\frac{\partial^2}{\partial\alpha^2}  l] &= -\beta \gamma\sum_{t=t_L}^{t_U} \frac{|t|^2}{(\alpha |t| +1)^{\beta+2}}
+ \gamma \beta(\beta+1) \sum_{t=t_L}^{t_U} \frac{|t|^2}{(\alpha |t| +1)^{\beta+2}}\\
&=\gamma \beta^2 \sum_{t=t_L}^{t_U} \frac{|t|^2}{(\alpha |t| +1)^{\beta+2}}.
\end{align*}
Note that $I_{\alpha\alpha}\rightarrow\infty$ as
$\gamma=E(y_0)\rightarrow\infty$.   However,  when $\gamma$ is fixed
\[
\lim_{\max(-t_L, t_U)\rightarrow\infty} I_{\alpha\alpha}=\infty \ \   \Leftrightarrow \ \ \beta \le 1
\ \   \Leftrightarrow \ \ \sum_{t=t_L}^{t_U} E\left( y_t\right) \rightarrow\infty.
\]

Next consider $I_{\alpha\beta}$.  Noting
\[
(\alpha t + 1)^{-\beta} = \exp(-\beta \log(\alpha t +1))
\]
and
\[
\frac{\partial}{\partial\beta}  (\alpha t + 1)^{-\beta} = - \log(\alpha t +1) \exp(-\beta \log(\alpha t +1))
= - \log(\alpha t +1) \frac{1}{(\alpha t + 1)^{\beta}},
\]
we have
\begin{align*}
-\frac{\partial^2}{\partial\alpha \partial\beta}  l
&= \frac{\partial}{\partial\alpha}
\left[
\sum_{t=t_L}^{t_U} y_t \log(\alpha |t| + 1)  - \gamma \sum_{t=t_L}^{t_U}  \log(\alpha |t| +1)\frac{1}{(\alpha |t| + 1)^{\beta}}
\right]\\
&= \sum_{t=t_L}^{t_U} y_t \frac{|t|}{\alpha |t| + 1} 
 - \gamma \sum_{t=t_L}^{t_U}  \frac{|t|}{(\alpha |t| + 1)^{\beta+1}}
+ \gamma\beta \sum_{t=t_L}^{t_U}\log(\alpha |t| + 1)  \frac{|t|}{(\alpha |t| + 1)^{\beta+1}}.
\end{align*}
When we take the expected value, the first two terms cancel and
\[
I_{\alpha\beta}=
E[-\frac{\partial^2}{\partial\alpha \partial\beta} l]=
\beta \gamma \sum_{t=t_L}^{t_U} \log(\alpha |t| + 1)  \frac{|t|}{(\alpha |t| + 1)^{\beta+1}}.
\]
The divergence is the same as in the case of $I_{\alpha\alpha}$.


Similarly we can evaluate $I_{\beta\beta}, I_{\alpha\gamma}, I_{\beta\gamma}, I_{\gamma\gamma}$ as
\begin{align*}
I_{\beta\beta}&=
\gamma \sum_{t=t_L}^{t_U} \frac{ (\log(\alpha |t| + 1))^2}{(\alpha |t| + 1)^\beta}, \qquad
  I_{\alpha\gamma}=-\beta\sum_{t=t_L}^{t_U} \frac{|t|}{(\alpha |t|+1)^{\beta+1}}, \\
I_{\beta\gamma} &= -\sum_{t=t_L}^{t_U}  \frac{\log(\alpha|t|+1)}{(\alpha |t| + 1)^{\beta}},\qquad 
I_{\gamma\gamma}= \frac{1}{\gamma}\sum_{t=t_L}^{t_U}  \frac{1}{(\alpha |t| + 1)^{\beta}}.
\end{align*}

Although it is tedious to prove the consistency and the asymptotic normality of
MLE based only on the computation of Fisher information matrix, our computation suggests
the following results.
Since $y_{t_0}$ has the Poisson distribution with the mean $\gamma$ and the standard
deviation $\sqrt{\gamma}$, $y_{t_0}/\gamma$ converges to 1 in probability as 
$\gamma\rightarrow\infty$. 
In fact, when we compute MLE, $\gamma$ is basically estimated by $y_{t_0}$, since the
number of postings has a sharp peak at $t=t_0$ (see Figure \ref{fig:sample1}). 
Furthermore $I_{\alpha\alpha}, I_{\alpha\beta},
I_{\beta\beta}$ are linear in $\gamma$.  This suggests that MLE is consistent 
as $\gamma\rightarrow\infty$.
For large expected value, Poisson distribution is approximately by 
normal distribution after normalization. Hence the score functions 
$\partial l/\partial\alpha,  \partial l/\partial\beta, \partial l/\partial\gamma$
are approximately normally distributed as $\gamma\rightarrow\infty$.
The confidence intervals given
in Section \ref{sec:data-analysis} are based on this approximation.
When $\gamma$ is fixed and $\max(-t_L, t_U)\rightarrow\infty$, the elements of the
Fisher information matrix diverge to $\infty$ if and only if $\beta\le 1$, or equivalently
$\sum_{t=t_L}^{t_U} E(y_t) \rightarrow\infty$.

\section{Conditional Poisson regression modeling for autocorrelations}
\label{sec:conditional-Poisson}

In the last section we assumed that the number of postings $y_t$ are independent.
We generalize this model to allow autocorrelations by conditional Poisson regression modeling. As we saw, the estimate of the parameter $\gamma$ in \eqref{eq:power-decay}
is very close to $y_{t_0}$.  Hence in this section we replace
$\gamma$ by $y_{t_0}$.  This is the initial value of our autoregressive scheme
and we model the number of postings after the event $y_t$, $t > t_0$.

Concerning the data $y_t$, $t<t_0$, before the date of the event, we can use
the model given in \eqref{eq:ar2} by reversing the time axis.
This is similar to look at
the standard AR(1) process
\[
x_t = \rho x_{t-1} + \epsilon_t
\]
in the reverse time direction by taking the
reciprocal of the autoregressive coefficient $\rho$. However this modeling of
the data before the date of the event is somewhat unsatisfactory, in particular
for the purpose of predicting $y_{t_0}$ before the date of the event.  We discuss this point 
again in Section \ref{sec:summary}.

We replace $\gamma$ by $y_{t_0}$ in \eqref{eq:power-decay} and
regard it 
as the conditional expected value of $y_t$ given $y_{t_0}$ 
\[
E(y_t | y_{t_0}) = y_{t_0}  \frac{1}{(\alpha |t-t_0|+1)^\beta}, \quad  t \ge  t_0.
\]
Then $E(y_t | y_{t_0})$ is recursively written as
\begin{align*}
E(y_t | y_{t_0}) &= E(y_{t-1} | y_{t_0})
  \left(\frac{(|t-t_0|-1)\alpha +1}{|t-t_0|\alpha +1}\right)^{\beta} \\
  &=E(y_{t-2}| y_{t_0})\left(\frac{(|t-t_0|-2)\alpha +1}{|t-t_0|\alpha +1}\right)^{\beta} \\
  &=  \dots
\end{align*}
We propose the following AR(2) type modeling of $y_t, t \ge t_0+2$:
\allowdisplaybreaks
\begin{align}
  p(y_{t_{0}+1}|y_{t_{0}})&={\rm Po}\left(y_{t_0+1} \mid y_{t_{0}}\frac{1}{(\alpha + 1)^\beta}\right) \nonumber \\
  p(y_{t_{0}+2}|y_{t_{0}+1},y_{t_{0}})&={\rm Po}\left(y_{t_0+2} \mid s\times y_{t_{0}+1}\left(\frac{\alpha + 1}{2\alpha + 1}\right)^\beta +
  (1-s)\times y_{t_{0}}\left(\frac{1}{2\alpha + 1}
  \right)^\beta \right) \nonumber
  \\
  p(y_{t_{0}+3}|y_{t_{0}+2},y_{t_{0}+1})&={\rm Po}\left(y_{t_0+3} \mid s\times y_{t_{0}+2}\left(\frac{2\alpha + 1}{3\alpha + 1}\right)^\beta +
  (1-s)\times y_{t_{0}+1}\left(\frac{\alpha+ 1}{3\alpha + 1}
  \right)^\beta \right) \nonumber
  \\
  &\ \vdots \nonumber
  \\
  p(y_{t_{0}+m}|y_{t_{0}+m-1},y_{t_{0}+m-2})&={\rm Po}\left(y_{t_0+m} \mid s\times y_{t_{0}+m-1}\left(\frac{(m-1)\alpha + 1}{m\alpha + 1}\right)^\beta  \right. \nonumber\\
  & \qquad + \left.
    (1-s)\times y_{t_{0}+m-2}\left(\frac{(m-2)\alpha+ 1}{m\alpha + 1}
  \right)^\beta \right). 
\label{eq:ar2}
\end{align}
Note that $y_{t_{0}+1}$ is given in an AR(1) form.  Then the conditional likelihood function for $\alpha,\beta,s$ given $y_{t_0}$ for the data $y_{t_0+1}, \dots, y_{t_0+T}$ is written as
\begin{align}
  L(\alpha,\beta,s) &=
  {\rm Po}\left(y_{t_0+1}\mid y_{t_{0}}\frac{1}{(\alpha + 1)^\beta}\right)
  \nonumber \\
  &\quad \times
  \prod_{m=2}^{T}
       {\rm Po}\left(y_{t_0+m} \mid s\times y_{t_{0}+m-1}\left(\frac{(m-1)\alpha + 1}{m\alpha + 1}\right)^\beta +
       (1-s)\times y_{t_{0}+m-2}\left(\frac{(m-2)\alpha+ 1}{m\alpha + 1}
       \right)^\beta \right).
\label{eq:ar2a}
\end{align}
When $s=1$, we have an AR(1) form.  When we estimate $s$ in $L(\alpha,\beta,s)$,
we restrict $s\in [0,1]$, although for some data sets unrestricted MLE of $s$
happened to be larger than 1.

Note that  the independence model in \eqref{eq:independent-poisson-likelihood}
and the AR(2) model in \eqref{eq:ar2a} are separate models.  In the usual
AR(1) model of continuous observations $x_t = \rho x_{t-1} + \epsilon_t$,
the independence model is a special case of $\rho=0$.
In order to interpolate between the independence model
\eqref{eq:independent-poisson-likelihood} and the AR(2) model
\eqref{eq:ar2a}, we propose the following more generalized and unifying model with 
the new parameters $u,v \in [0,1]$ representing the weights of the two models.
\begin{align}
&  p(y_{t_{0}+m}|y_{t_{0}+m-1},y_{t_{0}+m-2}) \nonumber \\
&\qquad = 
    {\rm Po}\Biggl(y_{t_{0}+m} \mid 
    w\times  \left( \frac{y_{t_{0}}}{((m-1)\alpha+1)^{\beta}}\right)^{u} \left(
    y_{t_{0}+m-1}\right)^{1-u}
    \times \left(\frac{(m-1)\alpha+1}{m\alpha+1}\right)^{\beta}
    \nonumber \\
&\qquad \ \ 
   +(1-w)\times
    \left( \frac{y_{t_{0}}}{((m-2)\alpha+1)^{\beta}}\right)^{v} \left(
    y_{t_{0}+m-2}\right)^{1-v}
    \times \left(\frac{(m-2)\alpha+1}{m\alpha+1}\right)^{\beta}
    \Biggl) .
\label{eq:war2}
\end{align}
In this unifying model, $u$ is the weight for the lag one term and $v$
is the weight for the lag two term.  We introduced these two parameters
separately for flexibility of the model.

\subsection{Fisher information matrix for the AR(1) model} 
Here we evaluate Fisher information matrix for AR(1) model, i.e.\ the model in \eqref{eq:ar2a} with $s=1$.  
We let $t_{0}=0$ for simplicity and assume that $y_0,\dots, y_T$ are observed.
We also replace $\gamma$ by $y_0$ and consider the conditional likelihood in $\alpha$
and $\beta$ given $y_0$.

The conditional expected values given $y_0$ are evaluated as
\begin{equation*}
  \begin{split}
    E[y_{1}|y_{0}]&=y_{0}\left(\frac{1}{\alpha +1}\right)^{\beta},
    \\
    E[y_{2}|y_{1}]&=y_{1}\left(\frac{\alpha+1}{2\alpha +1}\right)^{\beta},\ \ E^{y_{1}}[E[y_{2}|y_{1}]|y_0]=E[y_{2}|y_0]=y_{0}\left(\frac{1}{2\alpha +1}\right)^{\beta},\\
    &\vdots \\
    E[y_{t}|y_{t-1}]&=y_{t-1}\left(\frac{(t-1)\alpha+1}{t\alpha +1}\right)^{\beta},\ \ E^{y_{t-1}}[E[y_{t}|y_{t-1}]|y_0]=E[y_{t}|y_0]=y_{0}\left(\frac{1}{t\alpha +1}\right)^{\beta}.\\
  \end{split}
\end{equation*}
The conditional likelihood function is
\begin{equation*}
  L(\alpha, \beta)=\prod_{t=1}^{T}\frac{\mu_{t}^{y_t}}{y_{t}!}e^{-\mu_{t}},\quad  \mu_{t}=y_{t-1}\left(\frac{(t-1)\alpha+1}{t\alpha+1}\right)^{\beta}.
\end{equation*}
Then the conditional log-likelihood function $l(\alpha, \beta)=\log L(\alpha,\beta)$
is written as
\begin{equation*}
  \begin{split}
    l(\alpha, \beta)&=\sum_{t=1}^{T}(y_{t}\log\mu_{t}-\mu_{t}-\log y_{t}!) \\
    &= \sum_{t=1}^{T}\left\{y_{t}\left(\log y_{t-1}+\beta \log\frac{(t-1)\alpha+1}{ t\alpha+1} \right)-y_{t-1}\left(\frac{(t-1)\alpha +1}{t\alpha +1} \right)^{\beta}-\log y_{t}! \right\}\\
    &= C+\beta \sum_{t=1}^{T} y_{t}\log\frac{(t-1)\alpha+1}{ t\alpha+1}  -\sum_{t=1}^{T}y_{t-1}\left(\frac{(t-1)\alpha +1}{t\alpha +1} \right)^{\beta},
  \end{split}
\end{equation*}
where $C$ does not depend on $\alpha, \beta$.
The first derivative and the second derivative with respect to $\alpha$ are evaluated as
\begin{align*}
  -\frac{\partial}{\partial \alpha}l(\alpha, \beta)&
  =-\beta\sum_{t=1}^{T}y_{t} \left\{\frac{t-1}{(t-1)\alpha+1}-\frac{t}{t\alpha +1} \right\}-\beta\sum_{t=1}^{T}y_{t-1}\frac{((t-1)\alpha +1)^{\beta-1}}{(t\alpha+1)^{\beta+1}},
\\
    -\frac{\partial^2}{\partial \alpha^2}l(\alpha, \beta)&=-\beta\sum_{t=1}^{T}y_{t}\left\{-\frac{(t-1)^2}{((t-1)\alpha +1)^2}+\frac{t^2}{(t\alpha +1)^2} \right\}\\
    &\quad\;-\beta\sum_{t=1}^{T}y_{t-1}\frac{((t-1)\alpha+1)^{\beta-2}}{(t\alpha +1)^{\beta+2}}(-2\alpha t^2 -2(1-\alpha)t+1-\beta).
\end{align*}
Taking the expected value we have
\allowdisplaybreaks
\begin{align*}
I_{\alpha\alpha}
       &=-\beta y_0 \sum_{t=1}^{T}\frac{1}{(t\alpha +1)^{\beta}}\left\{-\frac{(t-1)^2}{((t-1)\alpha +1)^2}+\frac{t^2}{(t\alpha +1)^2} \right\}\\
    &\quad\;-\beta y_{0}\sum_{t=1}^{T}\frac{((t-1)\alpha+1)^{-2}}{(t\alpha +1)^{\beta+2}}(-2\alpha t^2 -2(1-\alpha)t+1-\beta)\\
    &=\beta y_{0}\sum_{t=1}^{T}\frac{((t-1)\alpha+1)^{-2}}{(t\alpha +1)^{\beta+2}}
    (-2\alpha t^2 -2(1-\alpha)t+1)\\
    &\quad\;+\beta y_{0}\sum_{t=1}^{T}\frac{((t-1)\alpha+1)^{-2}}{(t\alpha +1)^{\beta+2}}(2\alpha t^2 +2(1-\alpha)t-(1-\beta))\\
    &=\beta^{2}y_{0}\sum_{t=1}^{T}\frac{((t-1)\alpha+1)^{-2}}{(t\alpha +1)^{\beta+2}}.
\end{align*}
The mixed derivative with respect to $\alpha$ and $\beta$ and its expected value
are  evaluated as
\begin{align*}
    -\frac{\partial^2}{\partial \beta\partial\alpha}l(\alpha, \beta)&=-\sum_{t=1}^{T}y_{t} \left\{\frac{t-1}{(t-1)\alpha+1}-\frac{t}{t\alpha +1} \right\}\\
    &\quad\;-\sum_{t=1}^{T}y_{t-1}\frac{((t-1)\alpha +1)^{\beta-1}}{(t\alpha+1)^{\beta+1}}\left\{1+\beta\log\frac{(t-1)\alpha+1}{t\alpha +1}\right\},\\
I_{\alpha \beta}
   &=-y_{0}\sum_{t=1}^{T}\frac{1}{(t\alpha +1)^{\beta}} \left\{\frac{t-1}{(t-1)\alpha+1}-\frac{t}{t\alpha +1} \right\}\\
    &\quad\;-y_{0}\sum_{t=1}^{T}\frac{((t-1)\alpha +1)^{-1}}{(t\alpha+1)^{\beta+1}}\left\{1+\beta\log\frac{(t-1)\alpha+1}{t\alpha +1}\right\} \\
    &=y_{0}\sum_{t=1}^{T}\frac{((t-1)\alpha+1)^{-1}}{(t\alpha +1)^{\beta+1}}\\
    &\quad\;-y_{0}\sum_{t=1}^{T}\frac{((t-1)\alpha +1)^{-1}}{(t\alpha+1)^{\beta+1}}\left\{1+\beta\log\frac{(t-1)\alpha+1}{t\alpha +1}\right\}\\
    &=-\beta y_{0}\sum_{t=1}^{T}\frac{((t-1)\alpha +1)^{-1}}{(t\alpha+1)^{\beta+1}}\log\frac{(t-1)\alpha+1}{t\alpha +1}.
\end{align*}
Similarly, the second derivative with respect to $\beta$ and its expected values are  evaluated as
\begin{align*}
  -\frac{\partial^2}{\partial \beta^2}l(\alpha , \beta)&=\sum_{t=1}^{T}y_{t-1}\left(\frac{(t-1)\alpha+1}{t\alpha +1}\right)^{\beta}\left(\log\frac{(t-1)\alpha+1}{t\alpha
    +1}\right)^2,\\
I_{\beta\beta}
  &=y_{0}\sum_{t=1}^{T}\frac{1}{(t\alpha +1)^{\beta}}\left(\log\frac{(t-1)\alpha+1}{t\alpha +1}\right)^2.
\end{align*}

Note that the elements of the Fisher information matrix are proportional to $y_0$.
Also the relevant series diverges if and only if
$\beta\le 1$.  This is the same as in the power-law decay independence model.

For more complicated models of this section, the evaluation of Fisher information matrix is difficult, mainly because we can not separate $y_{t-1}$ and $y_{t-2}$ in
$\log \mu_t$.

\section{Data analysis of some Japanese social networking data}
\label{sec:data-analysis}
We apply our models to some social networking service data in Japan. The data which we used are summarized in Table \ref{tab:tab1}.  In Table \ref{tab:tab1}, ``Date'' is the date of the event in the format month/data in 2014.
``ID'' is our identifier of the events 
used in later tables.  ``Searchword'' is the word we used in BuzzFinder service to search for the postings related to the events.
``Remarks'' are the explanations of the events.

\subsection{Parameter estimation of the power-law decay independence model}
In Table 2 we show fits of the power-law decay independence model to data. 
Because many events showed asymmetry before and after the date of the event,
we estimated the before-event parameters $\alpha_{b},\beta_{b}$ and the peak level $\gamma$ 
for one week before the event, and then estimated the after-event  parameters $\alpha_{a},\beta_{a}$ 
separately with the same $\gamma$ as the before-event parameter.
We also computed  $95\%$ confidence intervals.

In the graph of Figure 3 we show the data of Valentine's day around February 14, 2014.
The graph looks almost symmetric at first sight, 
but the estimated before-event parameters
and after-event parameters were different. 
Indeed in the graph of Figure 3 the slope just before the date of the event is steeper than after the event
and the number of postings decrease to zero faster after the event than before the event.
Our estimated parameters reflect these facts.
\begin{figure}[htbp]
\centering
\begin{minipage}{0.45\columnwidth}
\includegraphics[width=6cm]{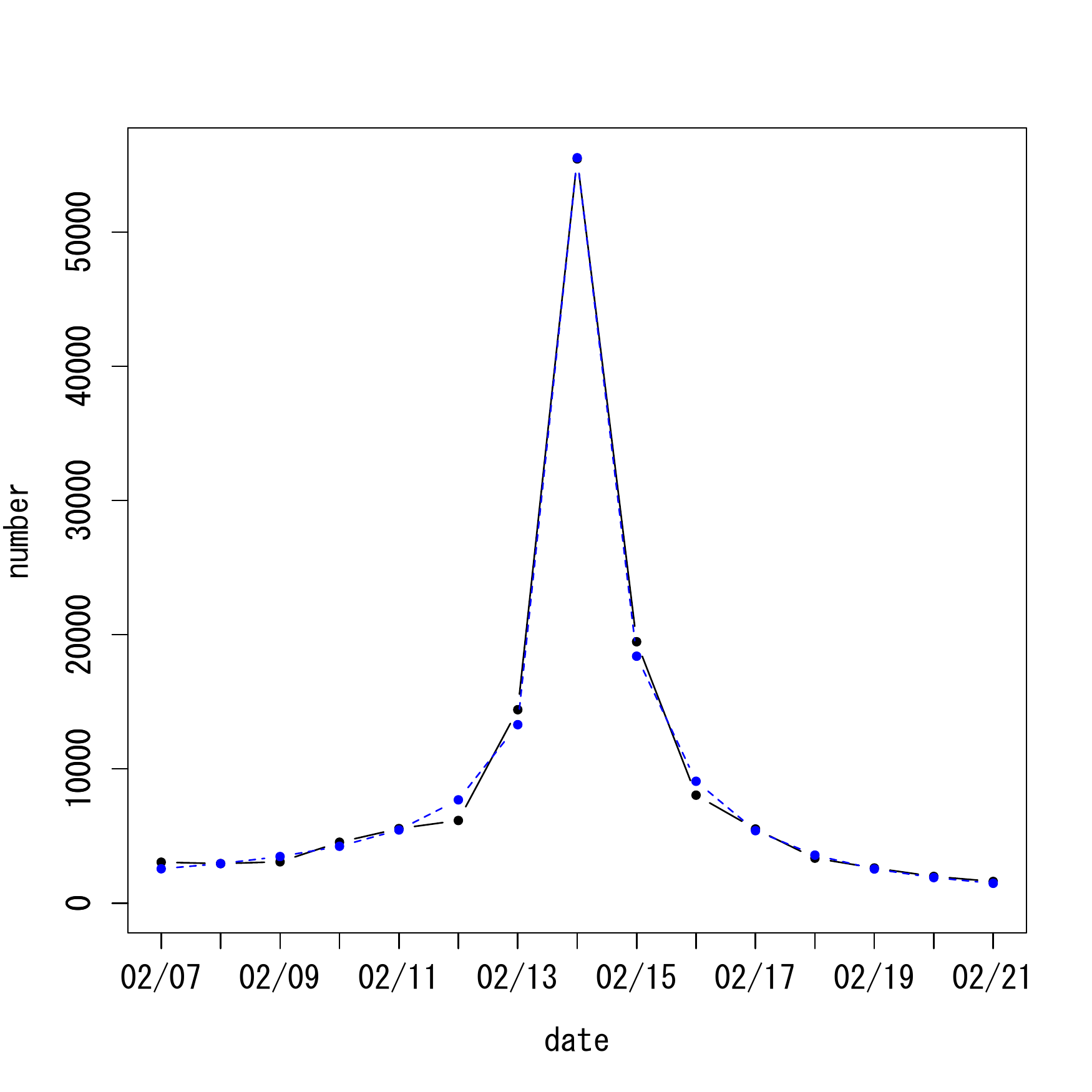}
\end{minipage}
\begin{minipage}{0.45\columnwidth}
\begin{tabular}{|c|c|}\hline
log-likelihood& $95\%$ confidence interval \\ \hline \hline
$\alpha_b=3.498$&$3.267<\alpha_b<3.728$ \\ \hline
$\beta_b=0.951$&$0.929<\beta_b<0.973$ \\ \hline
$\gamma=55538.8$&$54728.4<\gamma<56349.3$ \\ \hline
$\alpha_a=0.742$&$0.707<\alpha_a<0.777$ \\ \hline
$\beta_a=1.991$&$1.935<\beta_a<2.046$ \\ \hline
\end{tabular}
\end{minipage}
\caption{Valentine's Day 2014 for power-law decay independence model}
\label{fig:valetine}
\end{figure}

Based on the power-law decay independence model we considered predicting the after-event parameters based on the data before the event. However this was difficult, because of the asymmetry found in many events.
To confirm this phenomenon we performed multiple regression analysis, where
the before-event parameters $\alpha_{b},\beta_{b},\gamma$ are explanatory variables and 
the after-event parameters $\alpha_{a},\beta_{a}$ are objective variables.
But we did not find significant correlation.
  
\begin{figure}[htbp]
\centering
\includegraphics[width=6cm]{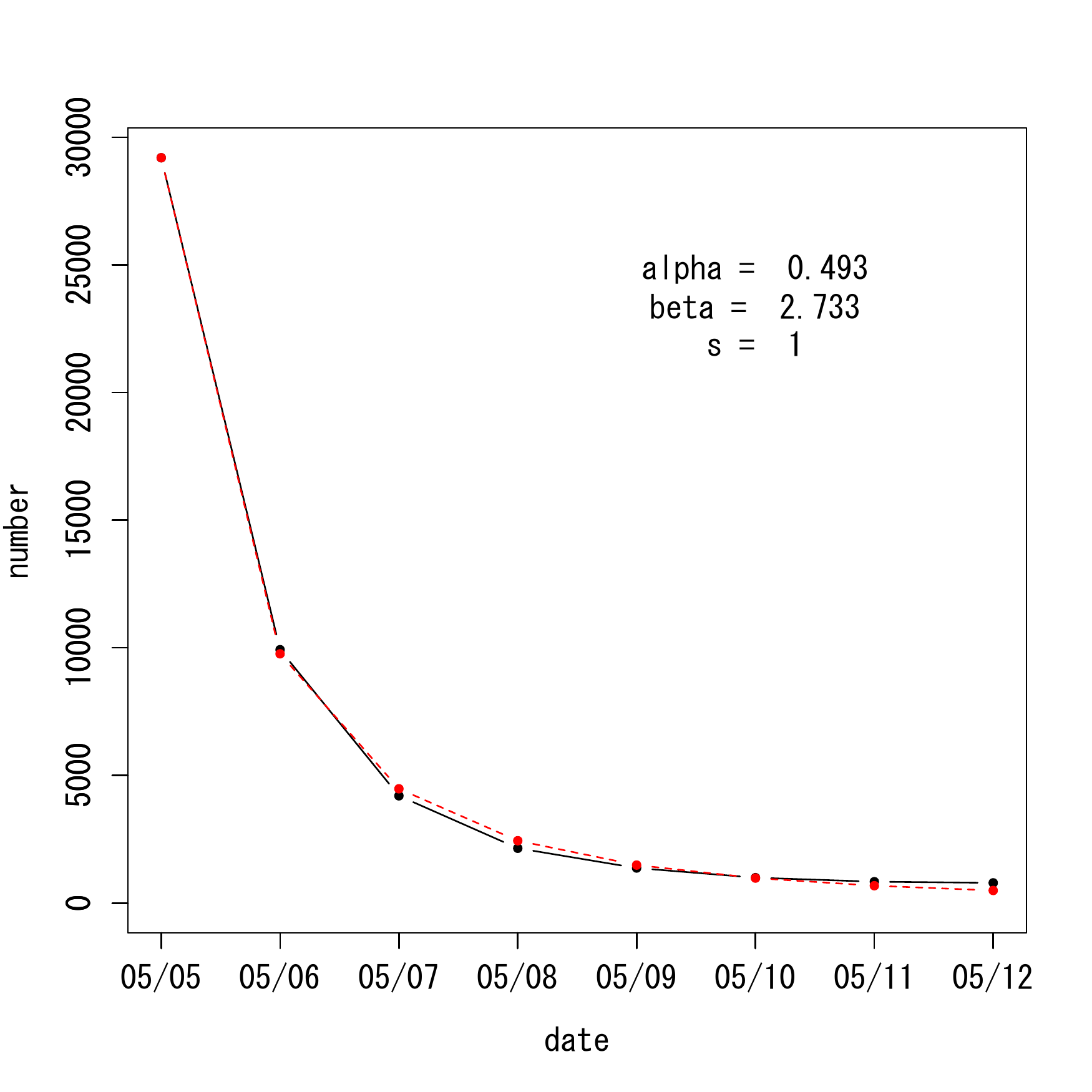} \quad
\includegraphics[width=6cm]{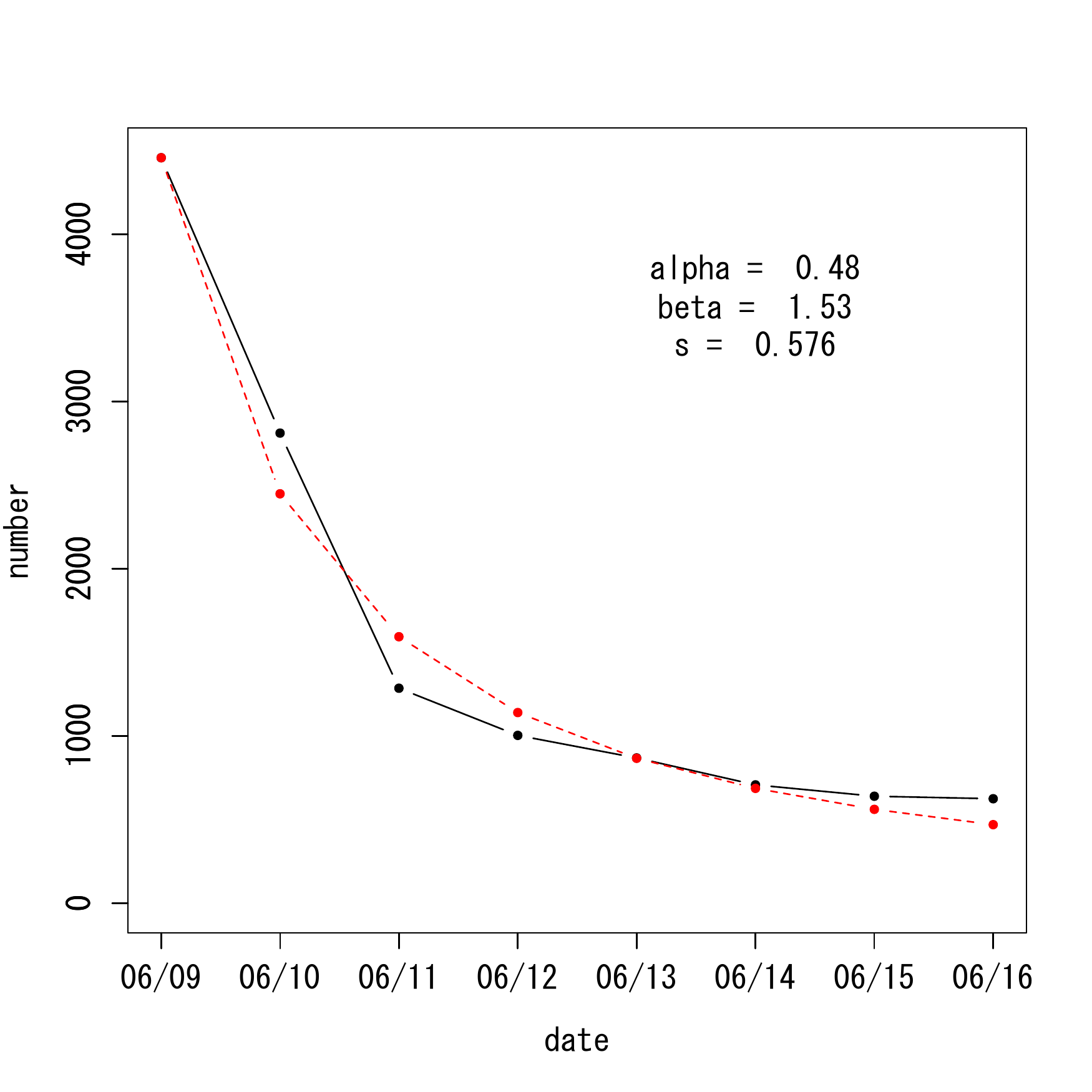}\\
Kodomo no Hi 2014 \hspace*{2.5cm} Yuko Oshima 
\caption{Fit of the AR(2) model}
\label{fig:ar2}
\end{figure}

\subsection{Parameter estimation of the AR(2) model}
In Table 3 we show the fit of AR(2) model to our data. In Table 3
``log-lik.'' stands for the log-likelihood of the estimated model. Figure 4 shows the fit of the AR(2) model
for ``Kodomo no Hi'' and for ``Yuko Oshima'' as typical examples.  
In Figure 4 $s$ is estimated as $s=1$ for ``Kodomo no Hi'', whereas $s$ is estimated
as $s=0.576$ for ``Yuko Oshima''.
It seems that the  parameter $s$ reflects the property of the event. The parameter
$s$ tends to be close to 1 for events with faster decay patterns, but tends to be less than one
for events with long-lasting interest after the events. 
This is reasonable, because $1-s$ represents the effect of two days ago 
and $s=1$ means that the autocorrelation is fully explained only by the number of 
postings one day ago. We compared AIC (Akaike's Information Criterion) 
for AR(2) model and AIC for AR(2) model with $s=1$. 
For many data sets AIC was smaller when $s$ is estimated to be less than 1.

\subsection{Parameter estimation of the unifying model}
In Table 4 we apply the unifying model \eqref{eq:war2} to data.
In the Table 5 we compare the unifying model and other models based on AIC. 
In many cases the values of the parameters $u, v$ are close to 1 in this model. 
This suggests that the unifying model is over-parameterized for many events and
the maximum likelihood estimation is not very stable.
Indeed when we compare AIC for various models, often other models have smaller AIC
than the unifying model.

\section{Summary and discussion}
\label{sec:summary}

In this paper we proposed a Poisson autoregression
model with the power-law decay of the mean parameter for the number of postings data.
Our model shows a good fit to various Japanese social networking data.
Also the parameters of our model are easy to interpret and our  model is
useful in describing patterns of the events.
Since the length of the data considered in this paper is fairly short,
covering only about one month, the unifying model in \eqref{eq:war2} with five parameters
is probably flexible enough. In fact for many events, we found that smaller models than
the unifying model showed better fits.
 
Our model assumes that there is a single date $t_0$ of an event.
Some events such as the Olympic games have longer duration. The pattern of
number of postings during the event with longer duration seems to be
more complicated, although the patterns before the beginning of the event 
and the after the event seem to be similar to single-day events.
It is not clear how to generalize our model to events with a
longer duration.

From practical viewpoint, it is important to predict the
peak level $\gamma$ of the number of postings and
the after-event parameters $\alpha_a, \beta_a$  before the event.
However we found that  this prediction was difficult for our data.  
Therefore in data analysis in Section \ref{sec:data-analysis}, we separately estimated
the before-event parameters and the after-event parameters, although for the modeling
purpose this is somewhat unsatisfactory.
As discussed at the beginning of Section \ref{sec:conditional-Poisson}
our conditional Poisson regression model is not suited for the prediction
before the event.  In addition, if the event has some surprising element on the
date of the event, then it is naturally difficult to predict it before the event.
We could use some characteristic of a particular event for the purpose of prediction.
For example, national holiday has a fixed date every year and we can analyze stability
of the pattern from year to year.


\begin{thebibliography}{8}
\bibitem{alfi-et-al}
Valentina Alfi, Andrea Gabrielli and Luciano Pietronero.
How people react to a deadline: time distribution of conference registrations and fee payments.
{\it Central European Journal of Physics}, {\bf 7}, No.3, 483--489, 2009.

\bibitem{1} 
Sitaram Asur, Bernardo A. Huberman, G\'abor Szab\'o  and Chunyan Wang.
Trends in Social Media : Persistence and Decay.  
In Proceedings of CoRR. pp.434--437.  2011. 

\bibitem{7} 
Patrick T. Brandt and John T. Williams.
A linear Poisson autoregressive model: the Poisson AR($p$) Model.
{\it Political Analysis}, {\bf 9}, Issue 2, 164--184, 2001.

\bibitem{fenner}
Trevor Fenner, Mark Levene and George Loizou.
A bi-logistic growth model for conference registration with an early bird deadline.
{\it Central European Journal of Physics}, {\bf 11}, No.7, 904--909, 2013.

\bibitem{fokianos-handbook}
Konstantinos Fokianos. Count time series models.
in {\it Time Series Analysis: Methods and Applications}, Handbook of Statistics, Vol.30, 
(ed. T.S.Rao, S.S.Rao and  C.R.Rao), North Holland, Amsterdam, 315--347, 2012.

\bibitem{fokianos-JASA} 
Konstantinos Fokianos,  Anders Rahbek,  and Dag Tj{\o}stheim.
Poisson autoregression. 
{\it Journal of the American Statistical Association}, {\bf 104}, 1430-1439, 2009.

\bibitem{fokianos-2011}
  Konstantinos Fokianos. Some recent progress in count time series.
{\it Statistics}, {\bf 45}, 49--58, 2011.

\bibitem{6} 
R.\ K.\ Freeland and B.\ P.\ M.\ McCabe. Analysis of low count time series data by Poisson autoregression, 
{\it Journal of Time Series Analysis},
{\bf 25}, Issue 5, 701--722, 2004.

\bibitem{3} Akshay Java, Xiaodan Song, Tim Finin and Belle Tseng. Why we twitter: understanding microblogging usage and communities,  in Proceedings of Joint 9th WEBKDD and 1st SNA-KDD Workshop '07, pp.56--65, 207.

\bibitem{5} 
Yukie Sano, Kenta Yamada, Hayafumi Watanabe, Hideki Takayasu and Misako Takayasu.
Empirical analysis of collective human behavior for extraordinary events in the blogosphere, {\it Physical Review E}, {\bf 87}, 
012805, 2013. 

\bibitem{4} Takeaki Sakaki, Makoto Okazaki and Yutaka Matsuo.  Earthquake shakes Twitter users: real-time event detection by social sensors. In  Proceedings WWW '10, pp.851--860. 2010.

\bibitem{traud}
Amanda L.\ Traud,  Peter J.\ Muchaa and Mason A.\ Porter.
Social structure of Facebook networks.
{\it Physica A}, {\bf 391}, 4165--4180, 2012. 

\bibitem{zhu} Fukang Zhu.
Modeling overdispersed or underdispersed count data with generalized
Poisson integer-valued GARCH models. {\it Journal of Mathematical Analysis and
Applications}, {\bf 389}, 58--71, 2012.
\end{thebibliography}

\begin{table}[thbp]
\begin{center}
\caption{Social networking data in Japan, 2014}
\label{tab:tab1}
\setlength{\tabcolsep}{3pt}
\begin{tabular}{|c|c|c|p{6cm}|}\hline
Date&ID&Searchword&Remarks \\ \hline \hline
1/13&Seijin&Seijin no Hi&Coming-of-Age Day \\ \hline
1/26&O-Marathon&Osaka International Ladies Marathon& \\ \hline
2/3&Setsubun&Setsubun&Bean Throwing Night \\ \hline
2/8&Sochi&Sochi Olympics&The opening ceremony took place February 7, local time. \\ \hline
2/9&Uemura&Aiko Uemura&The women's moguls final \\ \hline
2/11&Kenkoku&Kenkokukinenbi&National Foundation Day \\ \hline
2/14&Valentine&Valentine's Day& \\ \hline
2/15&Hanyu&Yuzuru Hanyu&Figure skating men's singles free \\ \hline
2/16&Kasai&Noriaki Kasai&Large hill individual men final \\ \hline
2/21&Asada&Mao Asada&Figure skating ladies's singles free \\ \hline
2/23&T-Marathon&Tokyo Marathon& \\ \hline
3/3&Academy&Academy Awards&The Oscar ceremonies took place March 2, local time. \\ \hline
3/11&Earthquake&The Great East Japan Earthquake& \\ \hline
3/14&White&White Day& \\ \hline
3/21&Shunbun&Shunbun no Hi&Vernal Equinox Day \\ \hline
3/24&mayoral&Osaka's mayoral elections& \\ \hline
5/5&Kodomo&Kodomo no Hi&Children's Day \\ \hline
6/9&Oshima&Yuko Oshima&AKB graduation performance \\ \hline
6/15&W-cup&The World Cup&In Japan's opening match against Cote d'Ivoire \\ \hline
6/21&Geshi&Geshi&the summer solstice \\ \hline
\end{tabular}
\end{center}
\end{table}

\begin{table}[htbp]
\begin{center}
\caption{Parameter estimation of the power-law decay independence model}
\begin{tabular}{|c|c|c|c|c|c|}\hline
ID&$\alpha_{b}$&$\beta_{b}$&$\gamma$&$\alpha_{a}$&$\beta_{a}$ \\ \hline \hline
Seijin&0.889&1.493&53664&0.908&1.463 \\ \hline
O-Marathon&7.902&0.858&993&0.207&4.337 \\ \hline
Setsubun&2.694&1.318&105334&0.397&2.921 \\ \hline
Sochi&0.179&5.499&5479&1.089&1.65 \\ \hline
Uemura&0.346&3.256&4992&0.231&3.436 \\ \hline
Kenkoku&1.405&2.012&7708&1.021&2.122 \\ \hline
Valentine&3.498&0.951&55539&0.742&1.991 \\ \hline
Hanyu&0.898&1.823&10317&2.353&0.891 \\ \hline
Kasai&599789.772&0.221&3368&0.131&4.584 \\ \hline
Asada&0.447&2.755&18898&1.289&1.047 \\ \hline
T-Marathon&2.327&1.186&8702&0.619&1.948 \\ \hline
Academy&39.489&0.423&2923&3.97&0.302 \\ \hline
Earthquake&10.695&0.822&41580&2.611&1.185 \\ \hline
White&6.285&0.691&45488&0.309&2.483 \\ \hline
Shunbun&2.99&1.259&14379&0.902&1.969 \\ \hline
Mayoral&0.096&7.803&1632&1.217&1.375 \\ \hline
Kodomo&14.37&0.783&29200&0.637&2.316 \\ \hline
Oshima&1.083&1.136&4456&0.673&1.23 \\ \hline
Wcup&0.013&34.707&45925&1.091&0.787 \\ \hline
Geshi&8.096&0.746&8183&0.766&1.498 \\ \hline
\end{tabular}
\end{center}
\end{table}

\begin{table}[htbp]
\begin{center}
\caption{Parameter estimation of the AR(2) model}
{\small
\begin{tabular}{|c|c|c|c|c|c|c|c|c|}\hline
&\multicolumn{6}{|c|}{AR(2) model}&\multicolumn{2}{|c|}{$s=1$} \\ \cline{2-9}
ID&$\alpha$&$\beta$&$\gamma$&$s$&log-lik.&AIC&log-lik.&AIC \\ \hline \hline
Seijin&0.612&1.834&53188&1.000&-342.303&690.606&-342.303&688.606 \\ \hline
O-Marathon&0.105&7.517&992&0.810&-58.919&123.839&-59.298&122.597 \\ \hline
Setsubun&0.265&3.944&105410&0.621&-2498.205&5002.410&-2617.576&5239.152 \\ \hline
Sochi&0.946&1.774&5374&0.671&-39.372&84.743&-39.575&83.150 \\ \hline
Uemura&0.276&3.173&5348&0.361&-89.743&185.485&-105.474&214.949 \\ \hline
Kenkoku&0.688&2.720&7679&1.000&-58.690&123.380&-58.690&121.380 \\ \hline
Valentine&0.722&2.040&55456&0.057&-154.633&315.266&-262.679&529.358 \\ \hline
Hanyu&1.919&0.951&9969&0.957&-82.052&170.104&-82.115&168.231 \\ \hline
Kasai&0.097&5.667&3381&0.227&-185.319&376.637&-246.634&497.269 \\ \hline
Asada&1.486&0.984&18939&0.341&-135.275&276.551&-181.722&367.444 \\ \hline
T-Marathon&0.331&2.934&8688&1.000&-185.955&377.911&-185.955&375.911 \\ \hline
Academy&0.692&0.743&2923&1.000&-958.221&1922.442&-958.221&1920.442 \\ \hline
Earthquake&2.615&1.207&41566&0.054&-298.935&603.871&-739.180&1482.360 \\ \hline
White&0.313&2.485&45438&0.131&-268.771&543.542&-459.159&922.317 \\ \hline
Shunbun&1.030&1.821&14369&0.405&-64.870&135.740&-73.371&150.742 \\ \hline
Mayoral&1.175&1.385&1581&0.342&-43.581&93.163&-49.383&102.766 \\ \hline
Kodomo&0.493&2.733&29201&1.000&-104.649&215.298&-104.649&213.298 \\ \hline
Oshima&0.480&1.530&4458&0.576&-155.904&317.808&-167.602&339.205 \\ \hline
Wcup&0.930&0.859&44480&0.666&-1994.083&3994.167&-2142.728&4289.457 \\ \hline
Geshi&0.552&1.789&8174&1.000&-140.557&287.115&-140.557&285.115 \\ \hline
\end{tabular}
}\end{center}
\end{table}

\begin{table}
\begin{center}
\caption{Parameter estimation of the unifying model}
\begin{tabular}{|c|c|c|c|c|c|c|c|c|}\hline
ID&$\alpha$&$\beta$&$\gamma$&$w$&$u$&$v$&log-lik.&AIC \\ \hline \hline
seijin&0.735&1.618&53188&0.643&0&1&-329.506&669.013 \\ \hline
O-Marathon&0.177&4.857&992&0.287&0.005&0.996&-53.204&116.407 \\ \hline
Setsubun&0.398&2.919&105410&0.998&0.999&0.999&-2092.256&4194.512 \\ \hline
Sochi&1.053&1.665&5374&1&1&1&-37.463&84.927 \\ \hline
Uemura&0.286&3.055&5348&1&1&1&-68.389&146.779 \\ \hline
Kenkoku&0.817&2.404&7679&0.670&0&1&-58.269&126.538 \\ \hline
Valentine&0.741&1.993&55456&0.007&0.988&0.561&-145.248&300.495 \\ \hline
Hanyu&2.168&0.901&9969&0.993&0.702&0.996&-71.639&153.279 \\ \hline
Kasai&0.135&4.488&3381&0&0.966&0.999&-138.561&287.121 \\ \hline
Asada&1.296&1.045&18939&0&1&1&-109.559&229.119 \\ \hline
T-Marathon&0.433&2.388&8688&0.674&0.002&0.995&-180.467&370.934 \\ \hline
Academy&4.524&0.288&2923&0.977&0.920&0.990&-604.371&1218.743 \\ \hline
Earthquake&2.710&1.173&41566&0.429&1&0&-277.531&565.061 \\ \hline
White&0.310&2.475&45438&0.010&0.997&0.873&-220.389&450.779 \\ \hline
Shunbun&0.901&1.969&14369&0&1&1&-59.885&129.770 \\ \hline
Mayoral&1.145&1.396&1581&0&1&1&-37.753&85.506 \\ \hline
Kodomo&0.493&2.733&29201&1&0&1&-104.649&219.298 \\ \hline
Oshima&0.674&1.229&4458&1&1&1&-120.807&251.614 \\ \hline
Ｗcup&0.966&0.815&44480&0&1&1&-1249.637&2509.273 \\ \hline
Geshi&0.632&1.619&8174&0.814&0&1&-137.969&285.939 \\ \hline
\end{tabular}
\end{center}
\end{table}

\begin{table}
\begin{center}
\caption{Comparison of models based on AIC}
\begin{tabular}{|c|c|c|c|c|c|}\hline
&&\multicolumn{4}{|c|}{$w=1$} \\ \cline{3-6}
&& \multicolumn{2}{|c|}{}&$u=1$&$u=0$ \\ \hline
ID&AIC&$u$&AIC&AIC&AIC \\ \hline \hline
Seijin&669.013&0.39&667.378&732.505&688.606 \\ \hline
O-Marathon&116.407&0.764&112.277&111.576&122.597 \\ \hline
Setsubun&4194.512&1&4190.483&4188.482&5239.152 \\ \hline
Sochi&84.927&1&80.927&78.927&83.15 \\ \hline
Uemura&146.779&1&142.777&140.777&214.949 \\ \hline
Kenkoku&126.538&0.334&122.515&127.391&121.38 \\ \hline
Valentaine&300.495&1&311.419&309.419&529.358 \\ \hline
Hanyu&153.279&0.705&149.28&150.843&168.231 \\ \hline
Kasai&287.121&1&283.097&281.095&497.269 \\ \hline
Asada&229.119&1&225.119&223.119&367.444 \\ \hline
T-Marathon&370.934&0.41&364.891&394.555&375.911 \\ \hline
Academy&1218.743&0.905&1212.032&1217.415&1920.442 \\ \hline
Earthquake&565.061&1&644.564&642.564&1482.36 \\ \hline
White&450.779&1&449.023&447.023&922.317 \\ \hline
Shunbun&129.77&1&125.77&123.77&150.742 \\ \hline
Mayoral&85.506&1&81.506&79.506&102.766 \\ \hline
Kodomo&219.298&0&215.298&273.882&213.298 \\ \hline
Oshima&251.614&1&247.613&245.613&339.205 \\ \hline
Ｗcup&2509.273&1&2505.273&2503.273&4289.457 \\ \hline
Geshi&285.939&0.211&282.618&357.08&285.115 \\ \hline
\end{tabular}
\end{center}
\end{table}

\end{document}